\documentclass[12pt,a4paper]{amsart}

\usepackage{amsmath,amssymb,amsthm}
\usepackage{graphicx}
\usepackage[dvips,ps,all]{xy}
\xyoption{line}
\xyoption{color}
\xyoption{frame}

\usepackage{wrapfig}
\usepackage{pxfonts}
\usepackage{multirow}

\newcommand{\LYYY}[5][{}]{\ensuremath{
 \xygraph{
!{<0pt,0pt>;<4pt,0pt>:<0pt,-4pt>::}
!{(1,4)}*{\scriptscriptstyle #1}
!{(1,3)}="a"
!{(1,1)}="b"
!{(-1,-1)}="c"
!{(3,-1)}="d"
!{(3,-2)}*{\scriptscriptstyle #5}
!{(-3,-3)}="e"
!{(1,-4)}*{\scriptscriptstyle #4}
!{(1,-3)}="f"
!{(-5,-5)}="g"
!{(-1,-5)}="h"
!{(-5,-6)}*{\scriptscriptstyle #2}
!{(-1,-6)}*{\scriptscriptstyle #3}
"b"-"c"
"b"-"d"
"c"-"e"
"c"-"f"
"e"-"g"
"e"-"h"
}
}}

\newcommand{\LYYRY}[5][{}]{\ensuremath{
 \xygraph{
!{<0pt,0pt>;<4pt,0pt>:<0pt,-4pt>::}
!{(1,4)}*{\scriptscriptstyle #1}
!{(1,3)}="a"
!{(1,1)}="b"
!{(-1,-1)}="c"
!{(3,-1)}="d"
!{(3,-2)}*{\scriptscriptstyle #5}
!{(-3,-3)}="e"
!{(-3,-4)}*{\scriptscriptstyle #2}
!{(1,-3)}="f"
!{(-1,-5)}="g"
!{(3,-5)}="h"
!{(-1,-6)}*{\scriptscriptstyle #3}
!{(3,-6)}*{\scriptscriptstyle #4}
"b"-"c"
"b"-"d"
"c"-"e"
"c"-"f"
"f"-"g"
"f"-"h"
}
}}

\newcommand{\LYRY}[5][{}]{\ensuremath{
 \xygraph{
!{<0pt,0pt>;<4pt,0pt>:<0pt,-4pt>::}
!{(1,4)}*{\scriptscriptstyle #1}
!{(1,3)}="a"
!{(1,1)}="b"
!{(-1,-1)}="c"
!{(3,-1)}="d"
!{(-2,-4)}*{\scriptscriptstyle #2}
!{(-2,-3)}="e"
!{(0,-4)}*{\scriptscriptstyle #3}
!{(0,-3)}="f"
!{(2,-3)}="g"
!{(4,-3)}="h"
!{(2,-4)}*{\scriptscriptstyle #4}
!{(4,-4)}*{\scriptscriptstyle #5}
"b"-"c"
"b"-"d"
"c"-"e"
"c"-"f"
"d"-"g"
"d"-"h"
}
}}

\newcommand{\RYYLY}[5][{}]{\ensuremath{
 \xygraph{
!{<0pt,0pt>;<4pt,0pt>:<0pt,-4pt>::}
!{(1,4)}*{\scriptscriptstyle #1}
!{(1,3)}="a"
!{(1,1)}="b"
!{(-1,-1)}="c"
!{(3,-1)}="d"
!{(5,-3)}="e"
!{(1,-3)}="f"
!{(-1,-5)}="g"
!{(3,-5)}="h"
!{(-1,-2)}*{\scriptscriptstyle #2}
!{(-1,-6)}*{\scriptscriptstyle #3}
!{(3,-6)}*{\scriptscriptstyle #4}
!{(5,-4)}*{\scriptscriptstyle #5}
"b"-"c"
"b"-"d"
"d"-"e"
"d"-"f"
"f"-"g"
"f"-"h"
}
}}

\newcommand{\RYYY}[5][{}]{\ensuremath{
 \xygraph{
!{<0pt,0pt>;<4pt,0pt>:<0pt,-4pt>::}
!{(1,4)}*{\scriptscriptstyle #1}
!{(1,3)}="a"
!{(1,1)}="b"
!{(-1,-1)}="c"
!{(3,-1)}="d"
!{(5,-3)}="e"
!{(1,-3)}="f"
!{(7,-5)}="g"
!{(3,-5)}="h"
!{(-1,-2)}*{\scriptscriptstyle #2}
!{(1,-4)}*{\scriptscriptstyle #3}
!{(3,-6)}*{\scriptscriptstyle #4}
!{(7,-6)}*{\scriptscriptstyle #5}
"b"-"c"
"b"-"d"
"d"-"e"
"d"-"f"
"e"-"g"
"e"-"h"
}
}}

\newcommand{\LYY}[4][{}]{\ensuremath{
 \xygraph{
!{<0pt,0pt>;<4pt,0pt>:<0pt,-4pt>::}
!{(1,4)}*{\scriptscriptstyle #1}
!{(1,3)}="a"
!{(1,1)}="b"
!{(-1,-1)}="c"
!{(3,-1)}="d"
!{(3,-2)}*{\scriptscriptstyle #4}
!{(-3,-3)}="e"
!{(-3,-4)}*{\scriptscriptstyle #2}
!{(1,-3)}="f"
!{(1,-4)}*{\scriptscriptstyle #3}
"b"-"c"
"b"-"d"
"c"-"e"
"c"-"f"
}
}}

\newcommand{\RYY}[4][{}]{\ensuremath{
 \xygraph{
!{<0pt,0pt>;<4pt,0pt>:<0pt,-4pt>::}
!{(5,-4)}*{\scriptscriptstyle #3}
!{(1,3)}="a"
!{(1,1)}="b"
!{(1,-4)}*{\scriptscriptstyle #2}
!{(-1,-1)}="c"
!{(3,-1)}="d"
!{(5,-3)}="e"
!{(-1,-2)}*{\scriptscriptstyle #4}
!{(1,-3)}="f"
!{(-1,-4)}*{\scriptscriptstyle #1}
"b"-"c"
"b"-"d"
"d"-"e"
"d"-"f"
}
}}

\newcommand{\LXX}[4][{}]{\ensuremath{
 \xygraph{
!{<0pt,0pt>;<4pt,0pt>:<0pt,-4pt>::}
!{(1,4)}*{\scriptscriptstyle #1}
!{(1,3)}="a"
!{(1,1)}*{\scriptscriptstyle{\bullet}}="b"
!{(2,1)}*{\scriptscriptstyle {}}
!{(-1,-1)}*{\scriptscriptstyle{\bullet}}="c"
!{(0,-1)}*{\scriptscriptstyle {}}
!{(3,-1)}="d"
!{(3,-2)}*{\scriptscriptstyle #4}
!{(-3,-3)}="e"
!{(-3,-4)}*{\scriptscriptstyle #2}
!{(1,-3)}="f"
!{(1,-4)}*{\scriptscriptstyle #3}
"b"-"c"
"b"-"d"
"c"-"e"
"c"-"f"
}
}}

\newcommand{\RXX}[4][{}]{\ensuremath{
 \xygraph{
!{<0pt,0pt>;<4pt,0pt>:<0pt,-4pt>::}
!{(5,-4)}*{\scriptscriptstyle #3}
!{(1,3)}="a"
!{(1,1)}*{\scriptscriptstyle{\bullet}}="b"
!{(1,-4)}*{\scriptscriptstyle #2}
!{(-1,-1)}="c"
!{(2,1)}*{\scriptscriptstyle {}}
!{(3,-1)}*{\scriptscriptstyle{\bullet}}="d"
!{(3,-2)}*{\scriptscriptstyle {}}
!{(5,-3)}="e"
!{(-1,-2)}*{\scriptscriptstyle #4}
!{(1,-3)}="f"
!{(-1,-4)}*{\scriptscriptstyle #1}
"b"-"c"
"b"-"d"
"d"-"e"
"d"-"f"
}
}}

\newcommand{\LZZ}[4][{}]{\ensuremath{
 \xygraph{
!{<0pt,0pt>;<4pt,0pt>:<0pt,-4pt>::}
!{(1,4)}*{\scriptscriptstyle #1}
!{(1,3)}="a"
!{(1,1)}*{\scriptscriptstyle{\circ}}="b"
!{(2,1)}*{\scriptscriptstyle {}}
!{(-1,-1)}*{\scriptscriptstyle{\circ}}="c"
!{(0,-1)}*{\scriptscriptstyle {}}
!{(3,-1)}="d"
!{(3,-2)}*{\scriptscriptstyle #4}
!{(-3,-3)}="e"
!{(-3,-4)}*{\scriptscriptstyle #2}
!{(1,-3)}="f"
!{(1,-4)}*{\scriptscriptstyle #3}
"b"-"c"
"b"-"d"
"c"-"e"
"c"-"f"
}
}}

\newcommand{\RZZ}[4][{}]{\ensuremath{
 \xygraph{
!{<0pt,0pt>;<4pt,0pt>:<0pt,-4pt>::}
!{(5,-4)}*{\scriptscriptstyle #3}
!{(1,3)}="a"
!{(1,1)}*{\scriptscriptstyle{\circ}}="b"
!{(1,-4)}*{\scriptscriptstyle #2}
!{(-1,-1)}="c"
!{(2,1)}*{\scriptscriptstyle {}}
!{(3,-1)}*{\scriptscriptstyle{\circ}}="d"
!{(3,-2)}*{\scriptscriptstyle {}}
!{(5,-3)}="e"
!{(-1,-2)}*{\scriptscriptstyle #4}
!{(1,-3)}="f"
!{(-1,-4)}*{\scriptscriptstyle #1}
"b"-"c"
"b"-"d"
"d"-"e"
"d"-"f"
}
}}

\newcommand{\LXZ}[4][{}]{\ensuremath{
 \xygraph{
!{<0pt,0pt>;<4pt,0pt>:<0pt,-4pt>::}
!{(1,4)}*{\scriptscriptstyle #1}
!{(1,3)}="a"
!{(1,1)}*{\scriptscriptstyle{\bullet}}="b"
!{(2,1)}*{\scriptscriptstyle {}}
!{(-1,-1)}*{\scriptscriptstyle{\circ}}="c"
!{(0,-1)}*{\scriptscriptstyle {}}
!{(3,-1)}="d"
!{(3,-2)}*{\scriptscriptstyle #4}
!{(-3,-3)}="e"
!{(-3,-4)}*{\scriptscriptstyle #2}
!{(1,-3)}="f"
!{(1,-4)}*{\scriptscriptstyle #3}
"b"-"c"
"b"-"d"
"c"-"e"
"c"-"f"
}
}}

\newcommand{\RXZ}[4][{}]{\ensuremath{
 \xygraph{
!{<0pt,0pt>;<4pt,0pt>:<0pt,-4pt>::}
!{(5,-4)}*{\scriptscriptstyle #3}
!{(1,3)}="a"
!{(1,1)}*{\scriptscriptstyle{\bullet}}="b"
!{(1,-4)}*{\scriptscriptstyle #2}
!{(-1,-1)}="c"
!{(2,1)}*{\scriptscriptstyle {}}
!{(3,-1)}*{\scriptscriptstyle{\circ}}="d"
!{(3,-2)}*{\scriptscriptstyle {}}
!{(5,-3)}="e"
!{(-1,-2)}*{\scriptscriptstyle #4}
!{(1,-3)}="f"
!{(-1,-4)}*{\scriptscriptstyle #1}
"b"-"c"
"b"-"d"
"d"-"e"
"d"-"f"
}
}}

\newcommand{\LZX}[4][{}]{\ensuremath{
 \xygraph{
!{<0pt,0pt>;<4pt,0pt>:<0pt,-4pt>::}
!{(1,4)}*{\scriptscriptstyle #1}
!{(1,3)}="a"
!{(1,1)}*{\scriptscriptstyle{\circ}}="b"
!{(2,1)}*{\scriptscriptstyle {}}
!{(-1,-1)}*{\scriptscriptstyle{\bullet}}="c"
!{(0,-1)}*{\scriptscriptstyle {}}
!{(3,-1)}="d"
!{(3,-2)}*{\scriptscriptstyle #4}
!{(-3,-3)}="e"
!{(-3,-4)}*{\scriptscriptstyle #2}
!{(1,-3)}="f"
!{(1,-4)}*{\scriptscriptstyle #3}
"b"-"c"
"b"-"d"
"c"-"e"
"c"-"f"
}
}}

\newcommand{\RZX}[4][{}]{\ensuremath{
 \xygraph{
!{<0pt,0pt>;<4pt,0pt>:<0pt,-4pt>::}
!{(5,-4)}*{\scriptscriptstyle #3}
!{(1,3)}="a"
!{(1,1)}*{\scriptscriptstyle{\circ}}="b"
!{(1,-4)}*{\scriptscriptstyle #2}
!{(-1,-1)}="c"
!{(2,1)}*{\scriptscriptstyle {}}
!{(3,-1)}*{\scriptscriptstyle{\bullet}}="d"
!{(3,-2)}*{\scriptscriptstyle {}}
!{(5,-3)}="e"
!{(-1,-2)}*{\scriptscriptstyle #4}
!{(1,-3)}="f"
!{(-1,-4)}*{\scriptscriptstyle #1}
"b"-"c"
"b"-"d"
"d"-"e"
"d"-"f"
}
}}

\newcommand{\lYsmall}[3]{\ensuremath{ 
 \xygraph{
!{<0pt,0pt>;<2.6pt,0pt>:<0pt,-2.6pt>::}
!{(-2,-3)}*{\scriptscriptstyle{#1}}
!{(2,-3)}*{\scriptscriptstyle{#2}}
!{(0,2)}="a"
!{(0,0)}="b"
!{(1.4,0)}*{\scriptscriptstyle{#3}}
!{(-2,-2)}="c"
!{(2,-2)}="d"
"b"-"c"
"b"-"d"
}
}}

\newcommand{\lY}[3]{\ensuremath{ 
 \xygraph{
!{<0pt,0pt>;<4.6pt,0pt>:<0pt,-4.6pt>::}
!{(-2,-3)}*{\scriptscriptstyle{#1}}
!{(2,-3)}*{\scriptscriptstyle{#2}}
!{(0,2)}="a"
!{(0,0)}="b"
!{(1.4,0)}*{\scriptscriptstyle{#3}}
!{(-2,-2)}="c"
!{(2,-2)}="d"
"b"-"c"
"b"-"d"
}
}}

\newcommand{\LYRYAB}[5][{}]{\ensuremath{
 \xygraph{
!{<0pt,0pt>;<5pt,0pt>:<0pt,-5pt>::}
!{(1,4)}*{\scriptscriptstyle #1}
!{(1,3)}="a"
!{(1,1)}*{\scriptscriptstyle{\circ}}="b"
!{(-1,-1)}*{\scriptscriptstyle{\bullet}}="c"
!{(3,-1)}*{\scriptscriptstyle{\circ}}="d"
!{(-2,-4)}*{\scriptscriptstyle #2}
!{(-2,-3)}="e"
!{(0,-4)}*{\scriptscriptstyle #3}
!{(0,-3)}="f"
!{(2,-3)}="g"
!{(4,-3)}="h"
!{(2,-4)}*{\scriptscriptstyle #4}
!{(4,-4)}*{\scriptscriptstyle #5}
"b"-"c"
"b"-"d"
"c"-"e"
"c"-"f"
"d"-"g"
"d"-"h"
}
}}

\DeclareMathAlphabet{\mathbbold}{U}{bbold}{m}{n}

\def\k{\mathbbold{k}}

\DeclareSymbolFont{rsfscript}{OMS}{rsfs}{m}{n}
\DeclareSymbolFontAlphabet{\mathrsfs}{rsfscript}

\DeclareFontFamily{OMS}{rsfs}{\skewchar\font'177}
\DeclareFontShape{OMS}{rsfs}{m}{n}{%
      <5> rsfs5
      <6> <7> rsfs7
      <8> <9> <10> rsfs10
      <10.95> <12> <14.4> <17.28> <20.74> <24.88> rsfs10
      }{}
\def\calF{\mathrsfs{F}}

\def\calK{\mathrsfs{K}}
\def\calL{\mathrsfs{L}}
\def\calM{\mathrsfs{M}}

\def\calP{\mathrsfs{P}}
\def\calQ{\mathrsfs{Q}}

\def\nn{\mathrm{no}\text{-}}

\DeclareMathOperator{\Ord}{Ord}
\DeclareMathOperator{\Vect}{Vect}

\DeclareMathOperator{\st}{st}

\makeatletter

\theoremstyle{plain}

\newtheorem {theorem}{Theorem}
\newtheorem{lemma}{Lemma}[theorem]
\newtheorem {corollary}{Corollary}
\newtheorem {conjecture}{Conjecture}

\newtheorem {proposition}{Proposition}

\theoremstyle{definition}

\newtheorem {definition}{Definition}
\newtheorem {remark}{Remark}
\newtheorem {example}{Example}

\begin{document}

\title{Pattern avoidance in labelled trees}
\author{Vladimir Dotsenko}
\address{Mathematics Research Unit, 
	 University of Luxembourg, 
	 Campus Kirchberg, 
	 6, Rue Richard Coudenhove-Kalergi, 
	 L-1359 Luxembourg, 
	 Grand Duchy of Luxembourg}
\email{vladimir.dotsenko@uni.lu}

\begin{abstract}
We discuss a new notion of pattern avoidance motivated by the operad theory: pattern avoidance in planar labelled trees. It is a generalisation of various types of consecutive pattern avoidance studied before: consecutive patterns in words, permutations, coloured permutations etc. The notion of Wilf equivalence for patterns in permutations admits a straightforward generalisation for (sets of) tree patterns; we describe classes for trees with small numbers of leaves, and give several bijections between trees avoiding pattern sets from the same class. We also explain a few general results for tree pattern avoidance, both for the exact and the asymptotic enumeration. 
\end{abstract}

\maketitle

\section{Introduction}

For pattern avoidance in words, apart from the ``real word interpretation'' (enumerate words not containing any obscene subwords), the pattern avoidance problem arises from studying (noncommutative) algebras with monomial relations. For example, describing words in the alphabet $\{\mathtt{A},\mathtt{B},\ldots,\mathtt{Z}\}$ not containing the word $\mathtt{FCUK}$ as a subword is equivalent to figuring out which monomials in generators $x_1,\ldots,x_{26}$ form a basis in the algebra $\k\langle x_1,\ldots,x_{26}\mid x_6x_3x_{21}x_{11}=0\rangle$. The significance of algebras with monomial relations is, in turn, explained by the theory of Gr\"obner bases which gives a method of finding a ``monomial replacement'' for every algebra with monomial relations~\cite{Ufn}.  Similarly, consecutive pattern avoidance in permutations~\cite{BS,Claesson} and coloured permutations~\cite{Ma} can be interpreted in terms of shuffle algebras with monomial relations~\cite{DKPerm}. The goal of this paper is to introduce to the combinatorics audience a new notion of pattern avoidance naturally arising when studying operads. Operads are similar to associative algebras, but while associative algebras and groups capture the kind of associativity that one observes when composing transformations of some set, operads capture the associativity exhibited when composing operations with several arguments. The property of an operation of having more than one argument results in a choice that is not present in the choice of algebras: one may either assume that our operations do not possess any symmetries or allow symmetries in the picture. In the former case, the pattern avoidance in question is the pattern avoidance in planar rooted trees; a few papers, both of combinatorial spirit \cite{Bacher,Parker,Rowland,GPPT} and operads-inspired~\cite{Loday} have been dealing with the arising questions of enumeration. In the case of operations with symmetries, the corresponding notion has not been studied before, and this paper is an attempt to give an elementary introduction to the arising research area. The ``right'' notion of a monomial relation for operations with symmetries is not as obvious as one might think: the action of symmetries makes every relation have too many consequences, and the arising class of ``operads with monomial relations'' appears to be way 
too narrow to be truly interesting or useful. The way to define monomial relations which avoids narrowing things down, which in particular led to a theory of operadic Gr\"obner bases, was suggested in~\cite{DK}; the corresponding algebraic object is called a shuffle operad. In this paper we, however, shall try to concentrate on the combinatorial aspects of the subject, touching the algebraic aspects only briefly; for details on the algebraic aspects, see, for example, \cite{DFree,DK1,DK}. Though we attempt to keep this article relatively self-contained, familiarity with the key results of the theory of consecutive pattern avoidance in permutations could be useful; for relevant historical information on pattern avoidance as well as the state-of-art of this area the reader is referred to a recent monograph~\cite{Kitaev}.

There are several types of questions in the theory of tree patterns which are meaningful from the operadic viewpoint. First of all, for a given set of patterns, exact enumeration results for trees avoiding that set are very important. Some examples of that sort appear in the following sections; in many cases the corresponding sequences of numbers are well known in combinatorics, but in many other cases one ends up with a sequences that appear to be unrelated to classical enumeration problems. Second, there is a question of asymptotic enumeration. We shall prove some results of that sort relying on the Golod--Shafarevich technique~\cite{GS} (which has recently been re-discovered in relation to combinatorics of pattern avoidance~\cite{GS1,GS2,GS3}). Also, there is a question of recognising the class of generating functions arising in this type of enumerating questions. However, while for word avoidance the answer is that the generating functions arising as answers for enumerating the avoidance of finite sets of words are always rational, even for consecutive pattern avoidance in permutations the class of arising generating functions does not have a satisfactory description. The answer is known if one ignores the leaf labels completely; in that case, the generating functions are algebraic (as proved in~\cite{Rowland} in the case of binary trees, in~\cite{GPPT} in the case of ternary trees, and in~\cite{KP} in the general case). For the notion of tree pattern avoidance discussed here, it follows from one of the theorems of~\cite{KP} that under some additional assumptions on the set of patterns (``shuffle regularity'') the generating functions for tree pattern avoidance are differentially algebraic (i.e. satisfy non-linear differential equations with polynomial coefficients). Finally, it is interesting to enumerate the Wilf equivalence classes of sets of tree patterns\footnote{Equivalence classes for ``strong equivalence'' might lead to interesting combinatorial results as well, but are much less natural from the operad point of view.}, and we discuss some basic results of that sort.

This paper is organised as follows. In Section~\ref{definitions}, we define planar labelled rooted trees and tree patterns, and show that this notion includes the classical types of consecutive pattern avoidance as particular cases. In Section~\ref{asymptotics}, we present some results on asymptotics for tree pattern avoidance, and introduce the notion of growth rate for a given set of patterns. In Section~\ref{clusters}, we formulate an exact enumeration result on tree pattern avoidance which follows from our work with Khoroshkin~\cite{DK1}, and discuss some consequences of that result and related questions. Finally, in Section~\ref{examples} we discuss various examples of patterns with small numbers of leaves. Each of these sections also contains some conjectures and natural questions that are beyond the scope of this paper. We also included an appendix explaining how pattern avoidance in trees arises in the operadic context. 

\subsection*{Acknowledgements. } I wish to thank Anton Khoroshkin and Dmitri Piontkovski for sending me a copy of their forthcoming preprint~\cite{KP}. 

\section{Planar labelled tree patterns}\label{definitions}

\subsection{Tree patterns}

\emph{Trees} have \emph{vertices} and \emph{edges}. A \emph{rooted tree} is a tree with a distinguished vertex, called the root. A rooted tree can be directed ``away from 
the root''; this way every vertex except for the root has exactly one \emph{parent}. Vertices whose common parent is a given vertex $v$ are called \emph{children} of $v$. 
Vertices that have at least one child are called \emph{internal}, vertices with no children are called \emph{leaves}. A \emph{planar rooted tree} is a rooted tree together 
with a total order on the set of children of each vertex, we shall think of it as of embedded in plane, the children of each vertex placed in the increasing order from the 
left to the right.

Throughout the paper, $X=\bigsqcup_{n\ge2} X_n$ is a finite alphabet represented as a disjoint union of its subsets $X_n$, $n\ge2$. 

\begin{definition}
A \emph{planar $X$-labelled rooted tree} is a planar rooted tree~$T$ with no internal vertices having exactly one child, and with a labelling of all vertices fulfilling the following restrictions: 
\begin{itemize}
 \item[-] every internal vertex~$v$ with $m$ children is labelled by an element $x_v\in X_m$;
 \item[-] every leaf of~$T$ is labelled by a positive integer in such a way that the following two conditions are satisfied:
 \begin{enumerate}
  \item (labelling set condition) the leaf labels are in one-to-one correspondence with the set $[l]=\{1,2,\ldots,l\}$ (where $l$ is the number of leaves of $T$);
  \item (local increasing condition) if we temporarily assign to each internal vertex $v$ the smallest of the labels of leaves that are descendants of~$v$ in~$T$  (thus every vertex of $T$, should it be an internal vertex or a leaf, has an integer assigned to it), then for each internal vertex $i$ the integers assigned to its children increase from the left to the right.
 \end{enumerate}
\end{itemize}
Notation: for $i\ge 0$, $l\ge 1$ we denote by $\mathcal{LT}_{i,l}(X)$ the set of planar $X$-labelled rooted trees with~$i$ internal vertices and $l$ leaves, 
\begin{equation}
\mathcal{LT}_{l}(X)=\bigcup_{i\ge 0}\mathcal{LT}_{i,l}(X) 
\end{equation}
is the set of planar $X$-labelled rooted trees with $l$ leaves, 
\begin{equation}
\mathcal{LT}(X)=\bigcup_{i\ge 0, l\ge 1}\mathcal{LT}_{i,l}(X) 
\end{equation}
is the set of all planar $X$-labelled rooted trees. In particular, $\mathcal{LT}_{0,l}(X)$ is non-empty only for $l=1$ (and in that case it consists of exactly one element, a one-vertex tree with no edges), and $\mathcal{LT}_{1}(X)=\mathcal{LT}_{0,1}(X)$ (since every internal vertex has at least two children, hence in the presence of an internal vertex we end up with at least two leaves).
\end{definition}

\begin{example}\label{ex:tree}
Let $X=X_2=\{\circ,\bullet\}$. The following tree is in $\mathcal{LT}_{7,8}(X)$:
\begin{equation}
\xygraph{
!{<0pt,0pt>;<7pt,0pt>:<0pt,-7pt>::}
!{(1,1)}*{\scriptscriptstyle{\bullet}}="b"
!{(-1,-1)}*{\scriptscriptstyle{\circ}}="c"
!{(3,-1)}="d"
!{(-3,-3)}*{\scriptscriptstyle{\bullet}}="e"
!{(1,-3)}*{\scriptscriptstyle{\circ}}="f"
!{(-4.4,-5)}*{\scriptscriptstyle{\circ}}="g"
!{(-1.9,-5)}="h"
!{(-0.4,-5)}*{\scriptscriptstyle{\circ}}="i"
!{(2.4,-5)}*{\scriptscriptstyle{\circ}}="j"
!{(-5.4,-7)}="k"
!{(-3.4,-7)}="l"
!{(-1.4,-7)}="m"
!{(0.4,-7)}="n"
!{(1.4,-7)}="o"
!{(3.4,-7)}="p"
!{(-5.4,-8)}*{\scriptscriptstyle 1}
!{(-3.4,-8)}*{\scriptscriptstyle 3}
!{(-1.4,-8)}*{\scriptscriptstyle 2}
!{(0.4,-8)}*{\scriptscriptstyle 7}
!{(1.4,-8)}*{\scriptscriptstyle 4}
!{(3.4,-8)}*{\scriptscriptstyle 6}
!{(-1.9,-6)}*{\scriptscriptstyle 5}
!{(3,-2)}*{\scriptscriptstyle 8}
"b"-"c"
"b"-"d"
"c"-"e"
"c"-"f"
"e"-"g"
"e"-"h"
"f"-"i"
"f"-"j"
"g"-"k"
"g"-"l"
"i"-"m"
"i"-"n"
"j"-"o"
"j"-"p"
} \ .
\end{equation}
\end{example}

A tree $T$ is said to be a \emph{left (right) comb} if for every internal vertex of~$T$ only its leftmost (rightmost) child may not be a leaf. Note that the local 
increasing condition makes the notions of ``left'' and ``right'' very different: for example, if a planar rooted tree~$t$ is a left comb, the only restriction on the leaf 
labels is that the leftmost leaf of $t$ is labelled by~$1$; by contrast, if $t$ is a right comb, there is only one leaf labelling satisfying the local increasing condition.

Most of the time throughout the paper we consider only the case $X=X_d$, thus assuming that for our trees all internal vertices have the same number of children (that is,~$d$ children). This assumption is mostly technical (it allows for closed formulas in various statements), in particular, all the asymptotic results we prove and conjecture in Section~\ref{asymptotics} are expected to be true in full generality. One particular simplification is that for such trees the number of internal vertices and the number of leaves are related: every tree with $k$ internal vertices has $kd-k+1$ leaves.

Let us prove a basic enumerative result which will be useful later. 

\begin{proposition}
For $X=X_d$, we have
\begin{equation}\label{number-of-trees}
|\mathcal{LT}_{kd-k+1}(X)|=|X_d|^k\frac{(kd)!}{(d!)^k\cdot k!}. 
\end{equation}
\end{proposition}

\begin{proof}
It would be useful for our purposes to consider, along with planar $X$-labelled rooted trees, two other types of trees. We denote by $\overline{\mathcal{T}}_n(X)$ the set of planar rooted trees with $n$ leaves whose internal vertices are labelled by~$X$, and by $\mathcal{T}_n(X)$ the set of planar rooted trees with $n$ leaves whose internal vertices are labelled by $X$, and whose leaves are labelled by $\{1,2,\ldots,n\}$ (in all possible ways). 

Recall that for $X=X_d=\{\bullet\}$, we have~\cite{Stanley2}
\begin{equation}
|\overline{\mathcal{T}}_{kd-k+1}(X)|=\frac{1}{kd-k+1}\binom{kd}{k}. 
\end{equation}
For the general case, we note that $|\overline{\mathcal{T}}_{kd-k+1}(X)|$ is $|X_d|^k$ times larger, since every of $k$ internal vertices should acquire a label from $X_d$. Also, it is clear that 
\begin{equation}
|\mathcal{LT}_{kd-k+1}(X)|=\frac{1}{(d!)^k}|\mathcal{T}_{kd-k+1}(X)|, 
\end{equation}
since for each of the $k$ internal vertices of a tree only one of the $d!$ permutations of its subtrees fulfils the local increasing condition. Finally, 
\begin{equation}
|\mathcal{T}_{kd-k+1}(X)|=(kd-k+1)!|\overline{\mathcal{T}}_{kd-k+1}(X)|, 
\end{equation}
since all leaf labelling are allowed in $\mathcal{T}_{kd-k+1}(X)$, and the formula follows.
\end{proof}

By a subtree of a planar labelled rooted tree~$T$ we always mean a subtree~$S$ with its root at one of the internal vertices of~$T$ such that for each internal vertex of $S$ all its children in $T$ are also its children in~$S$. (This way we can guarantee that the labels of internal vertices make sense for~$S$.)  Note that the second condition on leaf-labelling assigns to each internal vertex of a tree~$T$ a temporary integer label, so that a subtree~$S$ of a tree~$T$  almost belongs to $\mathcal{LT}(X)$: its leaves are labelled by integers such that the local increasing condition is satisfied (but the labelling set condition might not be satisfied). Replacing, for a subtree~$S$ with~$l$ leaves, its leaf labels by $1,2,\ldots,l$ in the unique order-preserving way, we shall obtain a tree $\st(S)\in\mathcal{LT}(X)$ which we call the standardisation of~$S$. 

\begin{definition}
A tree $T\in\mathcal{LT}(X)$ is said to \emph{contain a tree $P\in\mathcal{LT}(X)$ as a pattern} if there exists a subtree~$S$ of $T$ for which $P=\st(S)$. Otherwise $T$ is said to \emph{avoid}~$S$.
\end{definition}

\begin{example}\label{ex:pattern}
Let us recall the tree from Example~\ref{ex:tree}, and consider its subtree represented by thick lines in the following figure:
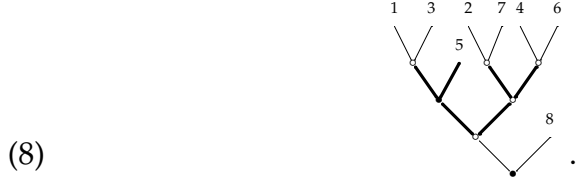
\begin{equation}
\xygraph{
!{<0pt,0pt>;<7pt,0pt>:<0pt,-7pt>::}
!{(1,1)}*{\scriptscriptstyle{\bullet}}="b"
!{(-1,-1)}*{\scriptscriptstyle{\circ}}="c"
!{(3,-1)}="d"
!{(-3,-3)}*{\scriptscriptstyle{\bullet}}="e"
!{(1,-3)}*{\scriptscriptstyle{\circ}}="f"
!{(-4.4,-5)}*{\scriptscriptstyle{\circ}}="g"
!{(-1.9,-5)}="h"
!{(-0.4,-5)}*{\scriptscriptstyle{\circ}}="i"
!{(2.4,-5)}*{\scriptscriptstyle{\circ}}="j"
!{(-5.4,-7)}="k"
!{(-3.4,-7)}="l"
!{(-1.4,-7)}="m"
!{(0.4,-7)}="n"
!{(1.4,-7)}="o"
!{(3.4,-7)}="p"
!{(-5.4,-8)}*{\scriptscriptstyle 1}
!{(-3.4,-8)}*{\scriptscriptstyle 3}
!{(-1.4,-8)}*{\scriptscriptstyle 2}
!{(0.4,-8)}*{\scriptscriptstyle 7}
!{(1.4,-8)}*{\scriptscriptstyle 4}
!{(3.4,-8)}*{\scriptscriptstyle 6}
!{(-1.9,-6)}*{\scriptscriptstyle 5}
!{(3,-2)}*{\scriptscriptstyle 8}
"b"-"c"
"b"-"d"
"c"-@[|(3)]"e"
"c"-@[|(3)]"f"
"e"-@[|(3)]"g"
"e"-@[|(3)]"h"
"f"-@[|(3)]"i"
"f"-@[|(3)]"j"
"g"-"k"
"g"-"l"
"i"-"m"
"i"-"n"
"j"-"o"
"j"-"p"
} \ .
\end{equation}
This subtree acquires the leaf numbering $\LYRYAB{1}{5}{2}{4}$, and after standardisation we get $\LYRYAB{1}{4}{2}{3}$. So $\LYRYAB{1}{4}{2}{3}$ occurs in our tree as a pattern. 
\end{example}

Throughout the paper, we only consider one type of tree patterns, so we often use the words ``tree pattern'' where one should say ``planar labelled rooted tree pattern''; 
we hope that the reader will appreciate this attempt to abbreviate things.

Let us fix some set of labels~$X$, and consider the pattern avoidance for planar $X$-labelled trees. The central question arising in the theory of pattern avoidance is that 
of enumeration of objects that avoid the given set of forbidden patterns~$\calP$ or, more generally, that contain 
exactly~$d$ occurrences of patterns from~$\calP$. This question naturally leads to the following equivalence relations for tree patterns. Two sets of tree patterns 
$\calP,\calP'\subset\mathcal{LT}(X)$ are said to be Wilf equivalent (notation: $\calP\sim_W\calP'$) if for every $l$, the number of $\calP$-avoiding trees with~$l$ leaves 
is equal to the number of $\calP'$-avoiding trees with~$l$ leaves. (In the case of (non-consecutive) permutation patterns the same notion was introduced by 
Wilf~\cite{Wilf}.) More generally, $\calP$ and $\calP'$ are said to be (strongly) equivalent (notation: $\calP\sim\calP'$) if for every $l$ and every $k\ge 0$, the number 
of trees with $l$ leaves that have $k$ occurrences of patterns from $\calP$ is equal to the number of trees with $l$ leaves that have $k$ occurrences of patterns from 
$\calP'$.

For enumeration, we shall primarily use the exponential generating functions with respect to the number of leaves in trees, so that, for example, the generating functions for the label set~$X$ and the pattern set $\calP$ are
$f_X(z)=\sum\limits_{n\ge 1}\frac{|X_n|z^n}{n!}$ and $f_\calP(z)=\sum\limits_{l\ge1}\frac{|\calP\cap\mathcal{LT}_{l}(X)|z^l}{l!}$ respectively. We shall denote the set of all trees avoiding the pattern set~$\calP$ by $\mathcal{LT}_{\nn{}\calP}(X)$, and its subset consisting of trees with $l$ leaves by $\mathcal{LT}_{l,\nn{}\calP}(X)$. The corresponding generating function is denoted by $f_{\nn{}\calP}(z)$:
\begin{equation}
f_{\nn{}\calP}(z)=\sum_{l\ge1}\frac{|\mathcal{LT}_{l,\nn{}\calP}(X)|z^l}{l!}. 
\end{equation}

\begin{remark}
One can also include the second variable in the generating series to count the internal vertices separately, and use the series
\begin{equation}
g_{\nn{}\calP}(z,t)=\sum_{i,l\ge1}\frac{|\mathcal{LT}_{i,l,\nn{}\calP}(X)|t^i z^l}{l!},
\end{equation}
which in particular would make all the sets finite when internal vertices of our trees are allowed to have a single child (and the label set $X$ includes~$X_1$). To keep the exposition simple, we avoid discussing these subtleties here. 
\end{remark}

The key feature of exponential generating functions in the context of planar $X$-labelled rooted trees is expressed by the following proposition.

\begin{proposition}\label{compos}
Suppose that $\calK$ and $\calL$ are two sets of planar $X$-labelled rooted trees. Let us define a set $\calM$ as follows: it consists of all trees $T$ that have an occurrence of a tree pattern from~$\calK$ rooted at the root of $T$, with the additional condition that all the subtrees rooted at the leaves of that pattern are occurrences of tree patterns from~$\calL$. Then
\begin{equation}
f_\calM(z)=f_\calK(f_\calL(z)).
\end{equation}
\end{proposition}

\begin{proof}
Clearly,
\begin{equation}
f_\calK(f_\calL(z))=\sum_{l\ge1}\frac{|\calK\cap\mathcal{LT}_{l}(X)|(f_\calL(z))^l}{l!},
\end{equation}
and it remains to note that the coefficient of $z^n$ in $(f_\calL(z))^l$ is the number of ordered forests of $l$ tree patterns from~$\calL$ with the total leaf set~$\{1,\ldots,n\}$, therefore $\frac{(f_\calL(z))^l}{l!}$ can be thought of as the enumerator for forests satisfying the increasing condition for minimal leaves. 
\end{proof}

\subsection{Tree patterns and other types of consecutive patterns}

In this section we assume, for simplicity, that $X=X_2$ (this corresponds to considering only binary tree patterns). 

Our first observation is that for $|X|=1$ the set of all permutations is naturally embedded in $\mathcal{LT}(X)$ as left combs: recall that a left comb has no conditions on where to put labels $2,3\ldots$, so left combs with $n+1$ leaves are in one-to one correspondence with permutations of length~$n$. If we denote by~$T(\sigma)$ the tree corresponding to the permutation~$\sigma$, then subtrees of $T(\sigma)$ are in one-to-one correspondence with subwords of $\sigma$, and the notion of a tree pattern for left combs corresponds precisely to the notion of a consecutive pattern for permutations. Moreover, if $\Pi$ is a set of consecutive permutation patterns, and $\calP_\Pi$ contains the left combs corresponding to elements of $\Pi$ and the right comb with three leaves, then the number of trees with $n+1$ leaves that avoid $\calP_\Pi$ is equal to the number of permutations of length~$n$ that avoid~$\Pi$. For $|X|>1$, the same construction with left combs leads naturally to the notion of pattern avoidance in coloured permutations~\cite{Ma}.

Moreover, the set of all words in a given alphabet $A$ is naturally embedded in $\mathcal{LT}(X)$ with $X=X_2=A$ as right combs. Indeed, recall that for a right comb with $n+1$ leaves there is exactly one way to label its leaves to fulfil the local increasing condition; to obtain a planar $X$-labelled rooted tree, it remains to label its internal vertices by~$A$, and the ways to do so are in one-to-one correspondence with $A$-words of length~$n$. If we denote by~$T(w)$ the tree corresponding to the word~$w$, then subtrees of $T(w)$ are in one-to-one correspondence with subwords of $w$, and the notion of a tree pattern for right combs corresponds precisely to the notion of a divisor for words. Moreover, if $W$ is a set of words, and $\calP_W$ consists of the right combs corresponding to elements of $W$ and all the left combs with three leaves, then the number of trees with $n+1$ leaves that avoid $\calP_W$ is equal to the number of words of length~$n$ without divisors from~$W$.

It is also possible to go the other way round and replace trees by objects resembling permutations and words. Let us assume, as above, that $|X|=1$ (ignoring this technical assumption will, as always, merely lead to coloured objects of the same sort). There is a very natural way to ``straighten'' the tree patterns and thus translate our questions into similar questions about patterns in sequences. Recall that the total number of planar labelled rooted tree patterns in our case is equal to $\frac{(2n)!}{2^n n!}=(2n-1)!!$, the double factorial number. This number is also equal~\cite{GeSt} to the number of permutations of 
the multiset $\{1,1,2,2,\ldots,n,n\}$ for which all the numbers appearing between the two occurrences of $k$ are greater than $k$ (for every $k=1,\ldots,n$). To a planar labelled rooted binary tree $T$, it is easy to assign recursively a permutation $\sigma(T)$ of that sort. For that, it is convenient to think of $T$ as of a left comb with subtrees grafted in the places of right children of internal vertices. We denote those subtrees by $T_1$, \ldots, $T_k$, in the order from the leftmost one to the rightmost one. Each subtree $T_i$ has its leaf labels belonging to a subset of $\{1,2,\ldots,n\}$, so, strictly speaking, they are not trees of the sort we consider, but, as usual, we can identify them with planar labelled rooted trees via standardisation, so we may apply $\sigma$ to them, obtaining permutations of appropriate multisets. We assume that $\sigma$ takes the only one-vertex tree to the empty word, and put
\begin{equation}
\sigma(T)=
\begin{cases}
1\sigma(T_1)1\quad \text{ for } k=1,\\ 
\st(\sigma(T_1)\sigma(T_2)\cdots\sigma(T_k))\quad \text{ for } k\ge 2.
\end{cases}
\end{equation}
For example, $\sigma(\lYsmall{1}{2}{})=11$, $\sigma(\LYRY{1}{2}{3}{4})=112332$, $\sigma(\LYYRY{1}{2}{3}{4})=122133$, $\sigma(\RYYLY{1}{2}{3}{4})=122331$; this leads to a meaningful notion of a generalised permutation pattern. It would be interesting to investigate this notion properly, in particular, to explore the links with patterns in set partitions \cite{JM,Kla,Sagan}, and also to see if the constructions of \cite{EKP} can be adapted here.

\section{Asymptotic enumeration}\label{asymptotics}

In this section, we discuss results on the asymptotic enumeration of trees avoiding a given set of patterns, where the results turn out to be in a way mimicking the results on the asymptotic enumeration for consecutive patterns in permutations~\cite{DKPerm,Eli1}. Our main tool is the following result, an adaptation of the classical Golod--Shafarevich inequality \cite{GS}; it is closely related to a similar inequality for symmetric operads~\cite{Piont}.

\begin{theorem}\label{GSh}
For every (possibly infinite) pattern set $\calP$, we have the following coefficient-wise inequality of power series:
\begin{equation}
f_\calP(f_{\nn{}\calP}(z))-f_X(f_{\nn{}\calP}(z))+f_{\nn{}\calP}(z)\ge z. 
\end{equation}
\end{theorem}

\begin{proof}
Let us consider two series of finite sets: the set $B_n$ is the subset of $\mathcal{LT}_{n}(X)$ consisting of trees $T$ whose subtrees rooted at the children of the root of~$T$ avoid patterns from~$\calP$, and the set $C_n$ is the subset of the set of pairs $\calP\times\mathcal{LT}_{n}(X)$ consisting of all pairs $(P,T)$ where there exists a subtree $S$ of $T$ rooted at the root of $T$ for which $\st(S)=P$, and all the trees rooted at the leaves of $S$ avoid patterns from~$\calP$. Proposition~\ref{compos} implies that
\begin{gather}
\sum_{n\ge 1}\frac{|B_n|}{n!}z^n=f_X(f_{\nn{}\calP}(z)), \\
\sum_{n\ge 1}\frac{|C_n|}{n!}z^n=f_\calP(f_{\nn{}\calP}(z)).
\end{gather}
Therefore, our power series inequality translates into 
\begin{equation}
|C_n|-|B_n|+|\mathcal{LT}_{n,\nn{}\calP}(X)|\ge 0, n\ge 2.
\end{equation}
This follows from an observation that there exists an obvious surjection from $C_n\sqcup \mathcal{LT}_{n,\nn{}\calP}(X)$ onto $B_n$: a tree from $B_n$ either avoids $\calP$, or has a pattern from $\calP$ rooted at its root. 
\end{proof}

\begin{corollary}\label{Lagrange}
Suppose that the power series $h(z)=\frac{1}{1-\frac{f_X(z)}{z}+\frac{f_\calP(z)}{z}}$ has non-negative coefficients. Then we have 
\begin{equation}
\frac{|\mathcal{LT}_{n,\nn{}\calP}(X)|}{n!} \ge \frac{1}{n}[z^{n-1}] h(z)^n.
\end{equation}
\end{corollary}

\begin{proof}
Let $g(z)=z-f_X(z)+f_\calP(z)$. According to the Lagrange inversion formula \cite{Stanley2}, for the coefficients of the compositional inverse $g^{\langle-1\rangle}(z)$ we have
\begin{equation}
[z^n] g^{\langle-1\rangle}(z)=\frac{1}{n} [z^{n-1}]\left(\frac{z}{g(z)}\right)^n= \frac{1}{n} [z^{n-1}] h(z)^n,
\end{equation}
so under our assumption on $h(z)$ the power series $g^{\langle-1\rangle}(z)$ has non-negative coefficients.
According to Theorem~\ref{GSh}, the power series 
\begin{equation}
f_\calP(f_{\nn{}\calP}(z))-f_X(f_{\nn{}\calP}(z))+f_{\nn{}\calP}(z)=g(f_{\nn{}\calP}(z)) 
\end{equation}
has non-negative coefficients as well, so we see that
\begin{multline}
f_{\nn{}\calP}(z)=g^{\langle-1\rangle}(g(f_{\nn{}\calP}(z)))=\\=
\sum_{n\ge1} \left(\frac{1}{n} [z^{n-1}] h(z)^n\right)g(f_{\nn{}\calP}(z))^n\ge
\sum_{n\ge1} \left(\frac{1}{n} [z^{n-1}] h(z)^n\right) z^n.
\end{multline}
This means that for every~$n\ge1$ we have
\begin{equation}
\frac{|\mathcal{LT}_{n,\nn{}\calP}(X)|}{n!}\ge \frac{1}{n} [z^{n-1}] h(z)^n,
\end{equation}
which is exactly what we want to prove.
\end{proof}

We use this result to prove the following theorem which gives a series of examples when the set of $\calP$-avoiding trees with~$n$ leaves grows rapidly, namely as $C^n$ times the number of all trees with $n$~leaves for some constant~$C$. A part of the proof is parallel to the corresponding proof in~\cite{DKPerm}.

\begin{theorem}
Suppose that $X=X_d$ for some~$d\ge2$, and that the power series 
\begin{equation}
h(z)=1-\frac{|X_d|z^{d-1}}{d!}+\frac{f_\calP(z)}{z}
\end{equation}
has a root $\alpha>0$. Then we have 
\begin{equation}
|\mathcal{LT}_{n,\nn{}\calP}(X)|\ge \left(\frac{(d!)^\frac{1}{d-1}}{|X_d|^\frac{1}{d-1}\alpha}\right)^{n-1} |\mathcal{LT}_{n}(X)|.
\end{equation}
\end{theorem}

\begin{proof}
Recall that if $X=X_d$, then the set $\mathcal{LT}_n(X)$ is non-empty only for $n=kd-k+1=k(d-1)+1$ for some~$k$. Therefore, $h(z)$ is a power series in $z^{d-1}$; we shall consider the series 
\begin{equation}
\hat{h}(t)=1-t+\sum_{t\ge 2}a_kt^k,
\end{equation}
for which $h(z)=\hat{h}\left(\frac{|X_d|}{d!}z^{d-1}\right)$; the only fact about the coefficients of that series that we use is that all the coefficients $a_k$ are non-negative ($a_k$ is a positive multiple of the number of certain labelled trees with~$k$ internal vertices). Under our assumption, this power series has a root $\beta=\frac{|X_d|}{d!}\alpha^{d-1}$. 

Let us consider multiplicative inverse, $\sum_{l\ge0}b_lt^l:=(\hat{h}(t))^{-1}$; clearly, $b_0=1$ and
$b_n-b_{n-1}+\sum_{k=2}^n a_kb_{n-k}=0$. Let us prove by induction that $b_n\ge \beta^{-1}b_{n-1}$. Indeed, for $n=1$ this statement is obvious ($\beta\ge1$ because otherwise $\hat{h}(\beta)$ is evidently positive), and for $n>1$ we note that by the induction hypothesis $b_{n-1}\ge\beta^{1-k}b_{n-k}$, so
\begin{multline*}
b_n=b_{n-1}-\sum_{k=2}^n a_kb_{n-k}\ge b_{n-1}-\sum_{k=2}^n a_k\beta^{k-1}b_{n-1}\ge\\ \ge
b_{n-1}-\sum_{k\ge2}a_k\beta^{k-1}b_{n-1}=\beta^{-1}b_{n-1}\left(\beta-\sum_{k\ge2}a_k\beta^k\right)=\beta^{-1}b_{n-1},
\end{multline*}
and the statement follows. Hence the $k^\text{th}$ coefficient of $(\hat{h}(t))^{-1}$ is greater than or equal to~$\beta^{-k}$. Therefore, the $(k(d-1))^\text{th}$ coefficient of~$h(z)$ is greater than or equal to
$\left(\frac{|X_d|}{d!}\right)\cdot\left(\frac{|X_d|\alpha^{d-1}}{d!}\right)^{-k}=\alpha^{-k(d-1)}$. 
This means that
\begin{equation}
h(z)\ge\sum_{k\ge1} \alpha^{-k(d-1)} z^{k(d-1)}=\left(1-\left(\frac{z}{\alpha}\right)^{d-1}\right)^{-1}, 
\end{equation}
so consequently
\begin{equation}
h(z)^n\ge\left(1-\left(\frac{z}{\alpha}\right)^{d-1}\right)^{-n}. 
\end{equation}
Since the coefficients of $h(z)$ are non-negative, Corollary \ref{Lagrange} applies, and we deduce that 
\begin{multline}
\frac{|\mathcal{LT}_{kd-k+1,\nn{}\calP}(X)|}{(kd-k+1)!}\ge \frac{1}{kd-k+1} [z^{k(d-1)}] h(z)^{kd-k+1} \ge\\ \ge
\frac{1}{kd-k+1}\alpha^{-k(d-1)}\binom{kd-k+1+k-1}{k}=\frac{1}{kd-k+1}\alpha^{-k(d-1)}\binom{kd}{k}.
\end{multline}
This, in the view of Formula~\eqref{number-of-trees}, simplifies to
\begin{multline}
|\mathcal{LT}_{kd-k+1,\nn{}\calP}(X)|\ge (kd-k+1)!\frac{1}{kd-k+1}\alpha^{-k(d-1)}\binom{kd}{k}=\\=
\alpha^{-k(d-1)} \frac{(kd)!}{k!}=
\left(\frac{d!}{\alpha^{d-1}}\right)^k \frac{(kd)!}{k!d!^k}=
\left(\frac{d!}{|X_d|\alpha^{d-1}}\right)^k |\mathcal{LT}_{kd-k+1}(X)|,
\end{multline}
which easily can be transformed into the inequality we want to prove.
\end{proof}

This theorem, in particular, can be used to obtain estimates in the case of one tree pattern in the case $|X|=1$, that is for trees with all internal vertices of the same arity and carrying the same label.

\begin{theorem}
Suppose that $|X|=1$, so that $X$ coincides with a one-element set $X_d$ for some~$d\ge2$, and that the set of forbidden patterns consists of one single pattern~$P$ with $k\ge 2$ internal vertices. Then for every pair $(d,k)$ except for $(2,2)$, $(2,3)$, $(2,4)$, and $(3,2)$ there exists a positive number $C$ (depending only on~ $d$ and~$k$ but not on the actual pattern) such that
\begin{equation}
|\mathcal{LT}_{n,\nn{}P}(X)|\ge C^{n-1} |\mathcal{LT}_{n}(X)|.
\end{equation}
\end{theorem}

\begin{proof}
It is enough to show that for all pairs $(d,k)$ except for $(2,2)$, $(2,3)$, $(2,4)$, and $(3,2)$, the polynomial $h(z)=1-\frac{z^{d-1}}{d!}+\frac{z^{k(d-1)}}{(kd-k+1)!}$ has a positive root. Denoting, as above, $t=\frac{z^{d-1}}{d!}$,  we see that it is enough to prove that the polynomial $\hat{h}(t)=1-t+\frac{d!^{k}t^k}{(kd-k+1)!}$ has a positive root. For that, we note that the derivative $-1+\frac{kd!^k t^{k-1}}{(kd-k+1)!}$ of the polynomial $\hat{h}(t)$ has the only positive root $t_0=\left(\frac{(kd-k+1)!}{k d!^k}\right)^\frac{1}{k-1}$; thus, to prove that $\hat{h}(t)$ has a positive root, it is enough to prove that its minimal value is attained at $t_0$ and is negative. Using the formula $t_0^{k-1}=\frac{(kd-k+1)!}{k d!^k}$, we see that $\hat{h}(t_0)=1-t_0+\frac{t_0}{k}$, so it suffices to prove that $t_0> \frac{k}{k-1}$, or
\begin{equation}
 \frac{(kd-k+1)!(k-1)^{k-1}}{d!^k k^k}> 1,
\end{equation}
which is true in all cases we consider, since it is true for $(d,k)=(2,5)$, $(d,k)=(3,3)$, and $(d,k)=(4,2)$, and its left hand side increases if either $d$ or~$k$ increases.
\end{proof}

Our results suggest a new numerical invariant of a set of patterns:

\begin{definition}
The \emph{growth rate} of a set of patterns $\calP\subset\mathcal{LT}(X)$ is 
\begin{equation}
\limsup_{n\to\infty} \left(\frac{|\mathcal{LT}_{n,\nn{}{}\calP}(X)|}{|\mathcal{LT}_{n}(X)|}\right)^{\frac{1}{n-1}}.
\end{equation}
\end{definition}

We expect that the methods of this section can be easily generalised to prove that for every set of labels $X$ (and for some $k$ depending on $X$) every pattern in $\mathcal{LT}(X)$ with at least $k$ internal vertices has positive growth. Moreover, computer experiments suggest that the following conjecture generalising Warlimont's conjecture for consectuive patterns~\cite{War} (proved recently by Ehrenborg, Kitaev and Perry \cite{EKP} via a beautiful link between consecutive pattern avoidance and the spectral theory for integral operators on the unit cube) holds for tree patterns. 

\begin{conjecture}
For every set of labels $X$, there exists an integer $d$ that for every pattern~$P$ in $\mathcal{LT}(X)$ with at least $d$ internal vertices and some numbers $c(P)>0$, $\lambda(P)>1$ we have
\begin{equation}
\frac{|\mathcal{LT}_{n,\nn{}P}(X)|}{|\mathcal{LT}_{n}(X)|}\sim c(P)\lambda(P)^{-n}.
\end{equation}
\end{conjecture}

It would be most interesting to adapt the approach of \cite{EKP} for tree patterns; such an adaptation, in addition to its consequences for the asymptotic enumeration questions, may be of substantial interest for the operad theory as well. Another natural question arising from the asymptotic enumeration is to understand the ``growth hierarchies'' of tree patterns, using the strategy of \cite{Kla1} or otherwise.

\section{Exact enumeration: cluster inversion formula}\label{clusters}

The following result on power series inversion is simultaneously a generalisation of the inversion formula for planar tree patterns \cite{Bacher,Loday,Parker} and of the cluster inversion formula of Goulden and Jackson for words and permutations~\cite{GJ,NZ}. This result is proved in \cite{DK1} by simple homological algebra; we formulate it here in a different way, so that an interested reader will easily prove it directly using the inclusion-exclusion formula, similarly to the usual cluster inversion.

\begin{definition}
Let $\calP$ be a set of patterns in $\mathcal{LT}(X)$. A tree $T$ together with a collection of its subtrees $T_1,\ldots,T_k$ is said to be a \emph{$k$-cluster} if the following two conditions hold:
\begin{enumerate}
 \item For every $i$, the subtree $T_i$ is an occurrence of a pattern from $\calP$: $\st(T_i)\in\calP$,
 \item for every edge $e$ of $T$ that joins two internal vertices, it is an edge between the two internal vertices of some $T_i$.
\end{enumerate} 
By definition, the set of $0$-clusters is the set of labels~$X$.
\end{definition}

Informally, a $k$-cluster is a tree which is completely covered by $k$ copies of patterns from $\calP$.

Let us denote by $c_{n,k}(\calP)$ the number of $k$-clusters for which the tree~$T$ has $n$ leaves.

\begin{theorem}[\cite{DK1}]
The compositional inverse for $f_{\nn{}\calP}(z)$ can be computed via clusters as follows:
\begin{equation}\label{eq:inverse}
f_{\nn{}\calP}^{\langle-1\rangle}(z)=z-\sum_{n\ge 1, k\ge 0}\frac{(-1)^k c_{n,k}(\calP)z^n}{n!}.
\end{equation}
\end{theorem}

The following consequence of the general inversion formula is a direct analogue for planar labelled tree patterns of the inversion formula \cite{Bacher,Loday,Parker} for planar tree patterns.

\begin{corollary}
Suppose that $X=X_2$, and that $\calP$ and $\calQ$ are two complementary sets of tree patterns with $3$ leaves, $\calP\sqcup\calQ=\mathcal{LT}_3(X)$. The corresponding generating functions for pattern avoidance are inverse to each other:
\begin{equation}\label{eq:inverse_q}
f_{\nn{}\calP}(-f_{\nn{}\calQ}(-z))=f_{\nn{}\calQ}(-f_{\nn{}\calP}(-z))=z.
\end{equation}
\end{corollary}

\begin{proof}
It is easy to see that in this case $k$-clusters for $\calP$ are in one-to-one correspondence with trees that avoid $\calQ$ (each such tree admits the unique covering by patterns from~$\calP$), and that the signs in the inversion formula match those suggested by~\eqref{eq:inverse_q} (since in our case the underlying tree of every $n$-cluster has $n+2$ leaves).
\end{proof}

Note that for left (right) combs corresponding to some consecutive permutation patterns (words), clusters for left combs are the left (right) combs corresponding to the usual Goulden--Jackson clusters for these patterns (words). This instantly proves the following result.

\begin{corollary}\label{trees-from-perm}
Suppose that two sets of consecutive permutation patterns (words) are Wilf equivalent. Then the two sets of tree patterns consisting of the left (right) combs corresponding to the given permutation patterns (words) are Wilf equivalent as tree patterns.
\end{corollary}

As in the case of permutations, we expect that at least in the case of a single pattern a careful study of its self-overlaps (``overlap sets'' of \cite{KS}, or equivalently ``overlap maps'' of~\cite{Nakamura}) would be very beneficial for studying Wilf equivalence. We shall discuss it in detail elsewhere, mentioning a particular case briefly in the next section.

Let us conclude this section with an open question. For consecutive patterns in permutations, Mendes and Remmel developed the symmetric functions method \cite{MR} for enumeration of permutations avoiding the given set of patterns. It is natural to expect that this method can be generalised to deal with the case of tree patterns, possibly making use of the plethysm for symmetric functions where our formulas compute compositions of power series. Some inversion formulas in the completion of the algebra of symmetric functions exist in the operadic context, being provided by homological algebra, in particular the operadic Koszul duality for symmetric operads~\cite{GK}, however, once we move from abstract trees to their representatives (that is, planar labelled rooted trees studied in this paper), there is no clear way to incorporate symmetric functions in the picture.

\section{Examples}\label{examples}

\subsection{Case $X=X_2$, $|X|=1$}

In this section, we assume that $X=X_2$, $|X|=1$, that is, we only work with binary trees, and do not use labels for internal vertices. To simplify the notation, we 
suppress $X$ in the formulas, and write simply $\mathcal{LT}_n$ etc.

We begin with classifying the Wilf classes of pattern sets with three leaves. 

\begin{theorem}\label{one-type}
There exist exactly four Wilf classes of sets of pattern sets with three leaves. 
\end{theorem}

\begin{proof}
This theorem follows from the two lemmas which give a precise description of the five Wilf classes. 
\begin{lemma}\label{one-pattern}
The three patterns $\LYY{1}{3}{2}{}{}$, $\LYY{1}{2}{3}{}{}$, $\RYY{2}{3}{1}{}{}$ are Wilf equivalent to each other; the number of trees with $n$ leaves avoiding either of them is equal to $(n-1)!$ for each $n\ge 3$.
\end{lemma}
\begin{proof}
Let us denote $P_1=\LYY{1}{3}{2}{}{}$, $P_2=\LYY{1}{2}{3}{}{}$, and $P_3=\RYY{2}{3}{1}{}{}$.

A correspondence $\rho$ between the set of $P_1$-avoiding trees and the set of $P_2$-avoiding trees can be defined recursively as follows. By definition, 
$\rho(\lYsmall{1}{2}{})=\lYsmall{1}{2}{}$. Let us represent a tree $T\in\mathcal{LT}_{n,\nn{}P_1}$ as a left comb with some subtrees $T_1,\ldots,T_k$ (listed along the way from the root) grafted at its ``right-looking'' leaves.  To determine $\rho(T)$, we apply $\rho$ to the subtrees $T_i$ and reverse the order of grafting. In other words, we graft $\rho(T_k)$ in the place of $T_1$, $\rho(T_{k-1})$ in the place of $T_2$ etc. Clearly, $\rho$ that identifies $\mathcal{LT}_{n,\nn{}P_1}$ with $\mathcal{LT}_{n,\nn{}P_2}$.

A correspondence $\kappa$ between the set of $P_1$-avoiding trees and the set of $P_3$-avoiding trees can be defined recursively as well. By definition, 
$\rho(\lYsmall{1}{2}{})=\lYsmall{1}{2}{}$. Let us represent a tree $T\in\mathcal{LT}_{n,\nn{}P_1}$ as a left comb with some subtrees $T_1,\ldots,T_k$ (listed along the way from the root) grafted at its right-looking leaves. We note that the set $\mathcal{LT}_{n,\nn{}P_3}$ is the set of all left combs with $n$ leaves; by induction we may assume that we already know the left combs $\kappa(T_1)$, \ldots, $\kappa(T_k)$. Let $\kappa(T)$ be the left comb whose right-looking leaves, listed along the way from the root, are the right-looking leaves of $\kappa(T_k)$, the right-looking leaves of $\kappa(T_{k-1})$, \ldots, the right-looking leaves of $\kappa(T_1)$. The observation (which can easily be proved by induction) that the label of the right-looking leaf of $\kappa(T)$ which is the farthest from the root is the equal to the smallest leaf label of $T_1$ shows how to construct the inverse of $\kappa$, so we identified $\mathcal{LT}_{n,\nn{}P_1}$ with $\mathcal{LT}_{n,\nn{}P_3}$.

In addition, since the set $\mathcal{LT}_{n,\nn{}P_3}$ is the set of all left combs with $n$ leaves, it has the cardinality~$(n-1)!$, so for each of the three subsets $\calP\subset\mathcal{LT}_{3}$ with $|\calP|=1$ and for each~$n$, there are exactly $(n-1)!$ different $\calP$-avoiding tree with $n$ leaves.
\end{proof}

\begin{lemma}
The three two-pattern sets $\{\LYY{1}{2}{3}{}{},\LYY{1}{3}{2}{}{}\}$, $\{\LYY{1}{2}{3}{}{},\RYY{2}{3}{1}{}{}\}$, and $\{\LYY{1}{3}{2}{}{},\RYY{2}{3}{1}{}{}\}$ are Wilf equivalent to each other; the number of trees with $n$ leaves avoiding either of these sets is equal to $1$ for each $n\ge 3$.
\end{lemma}
\begin{proof}
Indeed, for $\calP=\{\LYY{1}{2}{3}{}{},\LYY{1}{3}{2}{}{}\}$ the only $\calP$-avoiding tree with $n$ leaves is the only right comb, for $\calP=\{\LYY{1}{2}{3}{}{},\RYY{2}{3}{1}{}{}\}$ the only $\calP$-avoiding tree with $n$ leaves is the only left comb with labels of right-looking leaves increasing along the way from the root, and for $\calP=\{\LYY{1}{3}{2}{}{},\RYY{2}{3}{1}{}{}\}$ the only $\calP$-avoiding tree with $n$ leaves is the only left comb with labels of right-looking leaves decreasing along the way from the root. Therefore for each of the three subsets $\calP\subset\mathcal{LT}_{3}$ with $|\calP|=2$ and for each~$n$, there is exactly one $\calP$-avoiding tree with $n$ leaves. 

Alternatively, one can apply the inversion formula \eqref{eq:inverse_q}: from Lemma \ref{one-pattern}, we conclude that the exponential generating function for every one-pattern set is equal to $-\log(1-z)$; computing its inverse and adjusting the signs instantly shows that the exponential generating function for every two-pattern set is $\exp(z)-1$, which is the exponential generating function for the sequence $1,1,1,\ldots$.
\end{proof}

Since the empty pattern set and the pattern set containing all trees with three leaves form their own Wilf classes, the theorem follows.
\end{proof}

For pattern sets with at least four leaves, we have only partial results. Note that there are $15$ patterns with four leaves: $\LYYY{1}{2}{3}{4}$, $\LYYY{1}{2}{4}{3}$, $\LYYY{1}{3}{2}{4}$, $\LYYY{1}{3}{4}{2}$, $\LYYY{1}{4}{2}{3}$, $\LYYY{1}{4}{3}{2}$, 
$\LYYRY{1}{2}{3}{4}$, $\LYYRY{1}{2}{4}{3}$, $\LYYRY{1}{3}{4}{2}$, $\LYRY{1}{2}{3}{4}$, $\LYRY{1}{3}{2}{4}$, $\LYRY{1}{4}{2}{3}$, $\RYYLY{1}{2}{3}{4}$, $\RYYLY{1}{2}{4}{3}$, and $\RYYY{1}{2}{3}{4}$.
Therefore the number of sets of tree patterns with $4$ leaves is equal to $2^{15}=32768$, so a complete classification of Wilf classes is already a very heavy task. We shall present a very simple result on classification on Wilf classes for sets consisting of a single pattern.

\begin{theorem}
There exist exactly five Wilf classes of sets of one pattern with four leaves. 
\end{theorem}

\begin{proof}
This theorem follows from a sequence of lemmas which give a precise description of the five Wilf classes. 
\begin{lemma}\label{class1}
The three patterns $\LYYY{1}{2}{3}{4}$, $\LYYY{1}{4}{3}{2}$, and $\RYYY{1}{2}{3}{4}$ are Wilf equivalent to each other.  
\end{lemma}

\begin{proof}
There is a bijective proof which is completely analogous to that of Lemma~\ref{one-pattern}; we leave it to the reader to fill in the details.
\end{proof}

\begin{lemma}\label{class2}
The six patterns $\LYRY{1}{2}{3}{4}$, $\LYRY{1}{3}{2}{4}$, $\LYRY{1}{4}{2}{3}$, $\LYYRY{1}{3}{4}{2}$, $\RYYLY{1}{2}{3}{4}$, and $\RYYLY{1}{2}{4}{3}$ are Wilf equivalent to each other.  
\end{lemma}

\begin{proof}
Let us denote $P_1=\LYRY{1}{2}{3}{4}$, $P_2=\LYRY{1}{3}{2}{4}$, $P_3=\LYRY{1}{4}{2}{3}$, $P_4=\LYYRY{1}{3}{4}{2}$, $P_5=\RYYLY{1}{2}{3}{4}$, and $P_6=\RYYLY{1}{2}{4}{3}$.

Let us show that $P_1\sim_W P_2$ by exhibiting a one-to-one correspondence $\alpha$ between the $P_1$-avoiding patterns and the $P_2$-avoiding ones. If a tree $T$ avoids both $P_1$ and $P_2$, we put $\alpha(T)=T$. If $T$ avoids $P_1$ but contains $P_2$, we may assume that there is an occurrence of $P_2$ at the root of $T$ (otherwise we find the internal vertices \emph{closest to the root} that are roots of occurrences of $P_2$, and apply $\alpha$ recursively at these vertices). Let 
$T=\LYRY{T_1}{T_2}{T_3}{T_4}$, and denote $S_i=\alpha(T_i)$, $i=1,\ldots,4$. We put $\alpha(T)=\LYRY{S_1}{S_3}{S_2}{S_4}$. We constructed a bijection between $P_1$-avoiding trees containing $P_2$ and $P_2$-avoiding trees containing $P_1$. The case of $P_1$ and $P_3$ is handled in a similar way. 

Let us show that $P_1\sim_W P_4$ by exhibiting a one-to-one correspondence $\beta$ between the $P_1$-avoiding patterns and the $P_4$-avoiding ones. If a tree $T$ avoids both $P_1$ and $P_4$, we put $\beta(T)=T$. If $T$ avoids $P_1$ but contains $P_4$, we may assume that there is an occurrence of $P_4$ at the root of $T$ (otherwise we find the internal vertices \emph{closest to the root} that are roots of occurrences of $P_4$, and apply $\beta$ recursively at these vertices). Let 
$T=\LYYRY{T_1}{T_3}{T_4}{T_2}$, and denote $S_1=\beta(\lYsmall{T_1}{T_2}{})$, $S_2=\beta(\lYsmall{T_3}{T_4}{})$. We put $\beta(T)=\lYsmall{S_1}{S_2}{}$. Note that the only vertex of this tree where an occurrence of $P_4$ can be rooted is the root. However, if there is an occurrence of $P_4$ there, it is easily seen to imply an occurrence of $P_1$ in $T$, a contradiction. We constructed a bijection between $P_1$-avoiding trees containing $P_2$ and $P_2$-avoiding trees containing $P_1$. 

The equivalence $P_5\sim_W P_6$ can be established inductively similarly to how it is done in Lemma~\ref{one-pattern}.

Finally, the easiest way to see the equivalence $P_1\sim_W P_5$ is via the inverse generating functions. Basically, $P_1$ and $P_5$ have the same structure of self-overlaps: there are two self-overlaps, one of which is ``rigid'' (only one labelling of leaves is consistent with the local increasing condition), and the other one admits three different leaf labellings. This allows for an inductively constructed bijections between the clusters that control the coefficients of the inverse series. We leave the details to the reader. 
\end{proof}

\begin{lemma}\label{class3}
The four patterns $\LYYY{1}{2}{4}{3}$, $\LYYY{1}{3}{2}{4}$, $\LYYY{1}{3}{4}{2}$, and $\LYYY{1}{4}{2}{3}$ are Wilf equivalent to each other.  
\end{lemma}

\begin{proof}
The corresponding permutation patterns are Wilf equivalent, so the Corollary~\ref{trees-from-perm} applies. 
\end{proof}

Combining the lemmas above with a somewhat lengthy computation showing that 
\begin{itemize}
 \item for the class described in Lemma~\ref{class1} the sequence counting the trees avoiding the patterns of that class begins with $1,1,3,14,91,756,7657$,
 \item for the class described in Lemma~\ref{class2} the sequence counting the trees avoiding the patterns of that class begins with $1,1,3,14,90,739,7392$,
 \item for the class described in Lemma~\ref{class3} the sequence counting the trees avoiding the patterns of that class begins with $1,1,3,14,90,738,7364$,
 \item for the pattern $\LYYRY{1}{2}{4}{3}$ the sequence counting the trees avoiding that pattern begins with $1,1,3,14,90,740,7420$,
 \item for the pattern $\LYYRY{1}{2}{3}{4}$ the sequence counting the trees avoiding that pattern begins with $1,1,3,14,90,737,7336$,
\end{itemize}
we conclude that there are exactly five Wilf classes.
\end{proof}

Of the five integer sequences we saw in the previous proof, only two seem to appear in the Online Encyclopedia of Integer Sequences \cite{njas}: the third one matches 
\texttt{A088789}, the sequence of coefficients in the compositional inverse of the power series $\frac{2x}{1+\exp(x)}$, and the first one matches \texttt{A183611}, which, 
if we take care of the difference in numbering, is described as the sequence of coefficients of the power series $f(z)$ satisfying the differential equation 
$f''(z)=f'(z)^2+zf'(z)^3$. The first of these descriptions is not too surprising, as the inversion formula \eqref{eq:inverse} suggests that the inverse series makes lots of 
sense combinatorially. The second description is related to the results of Khoroshkin and Piontkovski \cite{KP} who proved that in some cases the generating function does 
indeed satisfy a differential equation; however, the patterns of that Wilf class are not covered by their results, so the appearance of a differential equation in this 
enumeration problem would be another bit of evidence supporting their general conjecture that states that for every finite set of patterns the corresponding generating function 
satisfies a non-linear differential equation with polynomial coefficients.

It is natural to ask which tree patterns are ``the hardest to avoid'', that have the fewest numbers of trees that avoid them, and which tree patterns are ``the easiest to avoid'', that is have the largest numbers of trees that avoid them. 
After examining the proof of the previous theorem and performing some computer experiments, we arrived at the following conjecture which is closely related to the conjecture of Elizalde and Noy~\cite{EN} that the permutation $12\ldots n$ (the identity permutation) is the easiest to avoid, and to the conjecture of Nakamura~\cite{Nakamura} that the permutation $12\ldots(n-2)n(n-1)$ (the transposition of its last two entries) is the hardest to avoid. 

\begin{conjecture}
Let us denote by $LC^<_n$, $LC^>_n$, and $RC_n$ the left comb with $n$ leaves whose leaf labels increase along the way from the root, the left comb with $n$ leaves whose leaf labels decrease along the way from the root, and the right comb with $n$ leaves respectively. The patterns $LC^<_n$, $LC^>_n$, and $RC_n$ are the easiest to avoid, and the pattern $\lY{RC_{n-1}}{n}{}$ is the hardest to avoid. 
\end{conjecture}

\subsection{Case $X=X_2$, $|X|=2$}

Throughout this section we assume that we work with binary trees with two possible labels for internal vertices, in other words, $X=X_2=\{\circ,\bullet\}$. 

The following result is still easy to obtain ``by hand''.

\begin{theorem}
There exist exactly two Wilf classes of sets of one pattern with three leaves, and exactly ten Wilf classes of sets of two patterns with three leaves. 
\end{theorem}

\begin{proof}
First of all, let us note that because of the inversion formula \eqref{eq:inverse_q}, we can work with the trees avoiding $11$ and $10$ patterns respectively. For the avoidance of $11$ patterns, the only allowed pattern can either use only one internal vertex label (in which case Theorem \ref{one-type} means that there exists only one allowed tree for each number of leaves) or use two different internal vertex labels (in which case there is no way to build an allowed tree with $4$ or more leaves). 

For the avoidance of $10$ patterns, our theorem follows from a sequence of lemmas which give a precise description of the ten Wilf classes. Lemmas \ref{easy1}--\ref{easy6} are almost obvious, since the trees that are being counted admit very explicit descriptions; we omit their proofs. Note that we switch to complements again: instead of listing the ten-element sets, we list their complements in the set of all tree patterns with three leaves.

\begin{lemma}\label{easy1}
The six sets whose complements in the set of all tree patterns with three leaves are $\{\LXX{1}{2}{3},\LXX{1}{3}{2}\}$, $\{\LXX{1}{3}{2},\RXX{2}{3}{1}\}$, $\{\LXX{1}{2}{3},\RXX{2}{3}{1}\}$, $\{\LZZ{1}{2}{3},\LZZ{1}{3}{2}\}$, 
$\{\LZZ{1}{3}{2},\RZZ{2}{3}{1}\}$, $\{\LZZ{1}{2}{3},\RZZ{2}{3}{1}\}$ are Wilf equivalent to each other; the number of trees with $n$ leaves avoiding either of these sets is equal to $(n-1)!$ for each $n\ge 4$.
\end{lemma}

\begin{lemma}\label{chainsaw}
The two sets whose complements in the set of all tree patterns with three leaves are $\{\LXZ{1}{2}{3},\RZX{2}{3}{1}\}$, $\{\LZX{1}{2}{3},\RXZ{2}{3}{1}\}$ are Wilf equivalent to each other; the number of trees with $n$ leaves avoiding either of these sets is equal to $\frac{3}{2}\left(\frac{n}{2}\right)!$ for even~$n\ge 4$ and to $\left(\frac{n-1}{2}\right)!+\frac12\left(\frac{n+1}{2}\right)!$ for odd $n\ge 3$.
\end{lemma}

\begin{lemma}\label{easy2}
The ten sets whose complements in the set of all tree patterns with three leaves are $\{\LXX{1}{2}{3}, \LXZ{1}{3}{2}\}$, $\{\LXX{1}{3}{2}, \LXZ{1}{2}{3}\}$, $\{\LXX{1}{2}{3}, \LZX{1}{3}{2}\}$, $\{\LXX{1}{3}{2}, \LZX{1}{2}{3}\}$, 
$\{\LZZ{1}{2}{3}, \LXZ{1}{3}{2}\}$, $\{\LZZ{1}{3}{2}, \LXZ{1}{2}{3}\}$, $\{\LZZ{1}{2}{3}, \LZX{1}{3}{2}\}$, $\{\LZZ{1}{3}{2}, \LZX{1}{2}{3}\}$, $\{\RXX{2}{3}{1}, \LZX{1}{2}{3}\}$, $\{\RZZ{2}{3}{1}, \LXZ{1}{2}{3}\}$ are Wilf equivalent to each other; the number of trees with $n$ leaves avoiding either of these sets is equal to $n-1$ for each $n\ge 4$.
\end{lemma}

\begin{lemma}\label{easy3}
The $32$ sets whose complements in the set of all tree patterns with three leaves are $\{\LXX{1}{2}{3}, \LZZ{1}{2}{3}\}$, $\{\LXX{1}{2}{3},\LZZ{1}{3}{2}\}$, $\{\LXX{1}{2}{3},\RZZ{2}{3}{1}\}$, $\{\LXX{1}{3}{2}, \LZZ{1}{2}{3}\}$,
$\{\LXX{1}{3}{2}, \LZZ{1}{3}{2}\}$, $\{\LXX{1}{3}{2}, \RZZ{2}{3}{1}\}$, $\{\RXX{2}{3}{1},\LZZ{1}{2}{3}\}$, $\{\RXX{2}{3}{1},\LZZ{1}{3}{2}\}$, $\{\RXX{2}{3}{1},\RZZ{2}{3}{1}\}$, $\{\LXX{1}{2}{3},\LZX{1}{2}{3}\}$, $\{\LXX{1}{3}{2},\LZX{1}{3}{2}\}$, $\{\LZZ{1}{2}{3},\LXZ{1}{2}{3}\}$, $\{\LZZ{1}{3}{2},\LXZ{1}{3}{2}\}$, $\{\RXX{2}{3}{1},\RZX{2}{3}{1}\}$, $\{\RZZ{2}{3}{1},\RZX{2}{3}{1}\}$, $\{\RXX{2}{3}{1},\RXZ{2}{3}{1}\}$, $\{\LXX{1}{3}{2},\LXZ{1}{3}{2}\}$, $\{\LZZ{1}{2}{3},\LZX{1}{2}{3}\}$, $\{\LZZ{1}{3}{2},\LZX{1}{3}{2}\}$, $\{\RXX{2}{3}{1},\RZX{2}{3}{1}\}$, $\{\RZZ{2}{3}{1},\RXZ{2}{3}{1}\}$, $\{\LXZ{1}{2}{3},\LZX{1}{2}{3}\}$, $\{\LXZ{1}{3}{2},\LZX{1}{3}{2}\}$, $\{\RXZ{2}{3}{1},\RZX{2}{3}{1}\}$, $\{\RXX{2}{3}{1}, \LZX{1}{3}{2}\}$, $\{\RZZ{2}{3}{1}, \LXZ{1}{3}{2}\}$, $\{\RZX{2}{3}{1}, \LXX{1}{2}{3}\}$, $\{\RZX{2}{3}{1}, \LXX{1}{3}{2}\}$, $\{\LZZ{1}{2}{3},\RXZ{2}{3}{1}\}$, $\LZZ{1}{3}{2},\RXZ{2}{3}{1}\}$, $\{\LXZ{1}{3}{2},\RZX{2}{3}{1}\}$, $\{\LZX{1}{3}{2},\RXZ{2}{3}{1}\}$
are Wilf equivalent to each other; the number of trees with $n$ leaves avoiding either of these sets is equal to $2$ for each $n\ge 4$,.
\end{lemma}

\begin{lemma}\label{easy4}
The two sets whose complements in the set of all tree patterns with three leaves are $\{\LZX{1}{2}{3},\RZX{2}{3}{1}\}$, $\{\LXZ{1}{2}{3},\RXZ{2}{3}{1}\}$
are Wilf equivalent to each other; the number of trees with $n$ leaves avoiding either of these sets is equal to~$1$ for $n=4$, and equal to~$0$ for each $n\ge 5$.
\end{lemma}

\begin{lemma}\label{easy5}
The two sets whose complements in the set of all tree patterns with three leaves are $\{\LZX{1}{3}{2},\RZX{2}{3}{1}\}$, $\{\LXZ{1}{3}{2},\RXZ{2}{3}{1}\}$ are Wilf equivalent to each other; the number of trees with $n$ leaves avoiding either of these sets is equal to $2$ for $n=4$, and equal to~$0$ for each $n\ge 5$.
\end{lemma}

\begin{lemma}\label{easy6}
The two sets whose complements in the set of all tree patterns with three leaves are $\{\LZX{1}{2}{3},\LZX{1}{3}{2}\}$, $\{\LXZ{1}{2}{3},\LXZ{1}{3}{2}\}$ are Wilf equivalent to each other; the number of trees with $n$ leaves avoiding either of these sets is equal to~$0$ for each $n\ge 4$.
\end{lemma}

\begin{lemma}\label{involutions}
The six sets whose complements in the set of all tree patterns with three leaves are $\{\LXX{1}{2}{3},\RXZ{2}{3}{1}\}$, $\{\LXX{1}{3}{2},\RXZ{2}{3}{1}\}$, $\{\LZZ{1}{2}{3},\RZX{2}{3}{1}\}$, $\{\LZZ{1}{3}{2},\RZX{2}{3}{1}\}$,
$\{\RXX{2}{3}{1},\LXZ{1}{3}{2}\}$, $\{\RZZ{2}{3}{1},\LZX{1}{3}{2}\}$ are Wilf equivalent to each other; the number of trees with $n$ leaves avoiding either of these sets is equal to the number of involutions in~$S_{n-1}$ for each $n\ge 4$.
\end{lemma}

\begin{proof}
Examining the structure of the allowed trees, we see that in each case we have a right comb or an increasing (decreasing) left comb, possibly with leaves replaced by pairs of leaves. Examining the leaf labels, one easily extract a decomposition of an involution into the product of disjoint cycles, which gives a one-to-one correspondence.
\end{proof}

\begin{lemma}\label{Fibonacci}
The two sets whose complements in the set of all tree patterns with three leaves are $\{\RXX{2}{3}{1},\LXZ{1}{2}{3}\}$ and $\{\RZZ{2}{3}{1},\LZX{1}{2}{3}\}$ are Wilf equivalent to each other; the number of trees with $n$ leaves avoiding either of these sets is equal to the $n^\text{th}$ Fibonacci number ($f_0=0$, $f_1=1$, $f_{n+1}=f_n+f_{n-1}$) for each $n\ge 3$.
\end{lemma}

\begin{proof}
Examining the structure of the allowed trees, we see that in each case we have a right comb, possibly with leaves replaced by pairs of leaves, but with the leaf labels increasing globally along the path from the root, therefore the number of trees we are trying to compute is equal to the number of sequences of $1$'s and $2$'s that sum up to~$n-1$, which is known to be the Fibonacci number~\cite{Stanley2}.
\end{proof}

\begin{lemma}  % $\frac{2(1+\sin(t))}{\cos(t)}-2-t$
The two sets whose complements in the set of all tree patterns with three leaves are $\{\LXZ{1}{3}{2},\LZX{1}{2}{3}\}$, $\{\LZX{1}{3}{2},\LXZ{1}{2}{3}\}$ are Wilf equivalent to each other; the number of trees with $n$ leaves avoiding either of these sets is twice the number of alternating permutations in $S_{n-1}$ for each $n\ge 3$.
\end{lemma}

\begin{proof}
Clearly, each allowed tree is a left comb, the labels $\circ$, $\bullet$ along the path from the root alternate, and so do the leaf labels. Therefore, the statement is obvious: alternating permutations come from the leaf labels, and ``twice'' reflect the fact that for each $k$ there are two alternating sequences of $\circ$'s and $\bullet$'s of length~$k$.
\end{proof}
\end{proof}

Further (computer-aided) investigation shows that the following statement is true:
\begin{theorem}
For sets of patterns with $3$ leaves, there are
\begin{enumerate}
 \item $2$ Wilf classes of $1$-element sets,
 \item $10$ Wilf classes of $2$-element sets,
 \item $40$ Wilf classes of $3$-element sets,
 \item $99$ Wilf classes of $4$-element sets,
 \item $189$ Wilf classes of $5$-element sets,
 \item $202$ Wilf classes of $6$-element sets,
 \item $189$ Wilf classes of $7$-element sets,
 \item $99$ Wilf classes of $8$-element sets,
 \item $40$ Wilf classes of $9$-element sets,
 \item $10$ Wilf classes of $10$-element sets,
 \item $2$ Wilf classes of $11$-element sets.
\end{enumerate}
\end{theorem}

To conclude this section, let us mention a promising direction towards new bijective proofs for enumeration of Wilf equivalence classes. In the case of unlabelled planar trees, many bijections have been constructed in~\cite{GPPT}, where the ``word notation'' for trees was used. That ``word notation'' is a construction of crucial importance for the operad theory. It was discovered by Hoffbeck~\cite{Hof} who defined a partial ordering compatible with the operad structure on the set of ``tree monomials'' which was then used to formulate and prove a PBW criterion for Koszul operads. Later, this partial ordering was extended to a total ordering (still compatible with the operad structure) in~\cite{DK}; that ordering is one of the key ingredients of the Gr\"obner bases machinery in the case of operads. The latter total ordering is coming from replacing a tree $T\in\mathcal{LT}(X)$ with a pair $(\mathsf{path}(T),\mathsf{perm}(T))$ consisting of a ``path sequence'' $\mathsf{path}(T)$ (a sequence of words in the alphabet~$X$) and a permutation $\mathsf{perm}(T)$ of the set of leaves of~$T$. All questions of pattern avoidance for trees can be translated into questions of pattern avoidance for this type of data; the corresponding notion of a pattern is a certain mixture of the classical notion of a divisor of a word and the notion of a generalised pattern in permutations~\cite{BS,Claesson} (different from the na\"ive notion of a generalised pattern in coloured permutations). This way of thinking about patterns in trees should be useful for bijective proofs; we hope to address it in more detail elsewhere. 

\appendix

\section{Shuffle operads and patterns in trees}

In this appendix, we recall some relevant definitions of the operad theory, and explain how the notion of pattern avoidance in trees arises naturally in this context.

Let us denote by~$\Ord$ the category whose objects are non-empty finite ordered sets (with order-preserving bijections as morphisms). Also, we denote by $\Vect$ the category of vector spaces (with linear operators as morphisms). It is usually enough to assume vector spaces finite-dimensional, though sometimes more generality is needed, and one assumes, for instance, that they are graded with finite-dimensional homogeneous components. 

\begin{definition}
\begin{enumerate}
\item A \emph{(non-symmetric) collection} is a contravariant functor from the category~$\Ord$ to the category~$\Vect$. We shall refer to images of individual sets as \emph{components} of our collection.
\item Let $\calP$ and $\calQ$ be two non-symmetric collections. The \emph{shuffle composition product} of $\calP$ and $\calQ$ is the non-symmetric collection $\calP\circ_{sh}\calQ$ defined by the formula
 $$
(\calP\circ_{sh}\calQ)(I):=\bigoplus_{k}\calP(k)\otimes\left(\bigoplus_{f\colon I\twoheadrightarrow[k]}\calQ(f^{-1}(1))\otimes\ldots\otimes\calQ(f^{-1}(k))\right),
 $$
where the sum is taken over all \emph{shuffling surjections}~$f$, that is surjections for which~$\min f^{-1}(i)<\min f^{-1}(j)$ whenever~$i<j$.
\item A shuffle operad is a monoid in the category of non-symmetric collections equipped with the shuffle composition product.
\end{enumerate}
\end{definition}

This definition is a counterpart of a more classical one, where one deals with finite sets without a specified order, and all surjections are allowed to define the composition. The corresponding monoids are called (symmetric) operads, and those are monoids widely used in algebra and topology, since they capture the algebraic properties of compositions of operations with several arguments. In particular, algebras of a certain type, like all associative algebras, or all Lie algebras, can be viewed as modules over a monoid of this kind (formed by all operations that can be defined on algebras of the given type), and this point of view proves to be useful. However, for computational purposes the fact that finite sets have symmetries gets in the way, and it turns out that shuffle operads allow to deal with some of the troubles arising because of that. (Every symmetric operad can be viewed as a shuffle operad, since every ordered set can be viewed as an unordered set.)

Monoids in the monoidal category of vector spaces (with the usual tensor product) are associative algebras. One can present associative algebras via generators and relations; if all relations are monomials in generators, the corresponding algebra admits a straightforward basis consisting of all monomials avoiding the monomials from the set of relations. In general, there is no such elegant description; to obtain it, one uses the machinery of Gr\"obner bases. A Gr\"obner basis is a special choice of a set of relations that allows to find a ``monomial replacement'' for the given algebra~$A$, that is an algebra with the same generators and with monomial relations for which the natural monomial basis is also a basis for~$A$. Such a set of monomial relations is provided by the ``leading terms'' of the relations forming a Gr\"obner basis.

It is very natural to try and find an appropriate Gr\"obner bases theory for operads. For operads with symmetries, it turns out to be impossible: not every symmetric operad has a monomial replacement. However, for shuffle operads it turns out to be possible. Let us be a little bit more precise. Similarly to the case of associative algebras, a shuffle operad can be presented via generators and relations, that is as a quotient of the free operad $\calF(V)$, where $V$ is the space of generators (which itself is a non-symmetric collection). If $X$ is the collection of ordered bases for components of $V$, that is a functor from $\Ord$ to $\Ord$, then the free shuffle operad generated by $V$ admits a basis of ``tree monomials'' which can be defined combinatorially; they are precisely planar $X$-labelled rooted trees studied throughout this paper. Any shuffle composition of tree monomials is again a tree monomial. The crucial feature of shuffle operads is that they admit good monomial orderings (and therefore one can talk about leading terms of relations). More precisely, there exist several ways to introduce a total ordering of tree monomials in such a way that all the shuffle compositions respect that total ordering: increasing one of the monomials we compose increases the result. Furthermore, the algebraic statement that a tree monomial~$T$ is obtained from another tree monomial~$S$ by shuffle compositions (that is, $T$ is divisible by~$S$ in our monoid) means, in the combinatorial language, that $S$ occurs in~$T$ as a pattern. 

We now have all the ingredients to relate questions about shuffle operads to pattern avoidance in trees. A Gr\"obner basis of an ideal $I$ in the free shuffle operad is a system $G$ of generators of~$I$ for which the leading monomial of every element of~$I$ is divisible by one of the leading terms of elements of~$G$. Such a system of generators allows to perform ``long division'' modulo~$I$, computing for every element its canonical representative. Tree monomials avoiding the leading terms of elements of~$G$ (``normal monomials'') form a basis in the quotient by the ideal~$I$; in other words, shuffle operads do admit monomial replacements. Thus, it is clear that enumerating trees avoiding the given set of patterns literally corresponds to computing dimensions of components of operads presented by generators and relations. (Of course, the first step is to compute a Gr\"obner basis, see~\cite{DK} for an algorithm approaching this problem.) From this point of view, the asymptotic enumeration is also very meaningful: in fact, the area of the operad theory primarily concerned with quotients of the operad describing associative algebras (better known as ``theory of varieties of associative algebras'') has asymptotic questions (``codimension growth'') among its core questions~\cite{GZ,KR}. There were even some attempts to use the Gr\"obner bases formalism in that context, see for example~\cite{Lat}. We expect that our approach will prove useful in that context. 

\bibliographystyle{amsplain}

\providecommand{\bysame}{\leavevmode\hbox to3em{\hrulefill}\thinspace}
\providecommand{\MR}{\relax\ifhmode\unskip\space\fi MR }
\providecommand{\MRhref}[2]{%
  \href{http://www.ams.org/mathscinet-getitem?mr=#1}{#2}
}
\providecommand{\href}[2]{#2}

\end{document}